\documentclass[12pt]{article}

\usepackage{array} 
\usepackage{amsfonts}
\usepackage{amsmath}
\usepackage{amssymb}
\usepackage{graphicx}
\usepackage{color}
\usepackage{setspace}
\newtheorem{thm}{Theorem}[section]

\newtheorem{lem}[thm]{Lemma}

\numberwithin{equation}{section}

\newcommand{\RR}{\mathbb{R}}

\newcommand{\ve}{\varepsilon}

\newcommand{\ren}{\RR^N}

\def\qed{\,\unskip\kern 6pt \penalty 500
\raise -2pt\hbox{\vrule \vbox to8pt{\hrule width 6pt
\vfill\hrule}\vrule}\par}
\definecolor{darkblue}{rgb}{0.05, .05, .65}
\definecolor{darkgreen}{rgb}{0.1, .65, .1}
\definecolor{darkred}{rgb}{0.8,0,0}

\begin{document}

\title{\textbf{The Dirichlet Problem for the fractional \\ p-Laplacian evolution equation}}

\author{\Large ~by~ Juan Luis V\'azquez
\footnote{Departamento de Matem\'{a}ticas, Universidad Aut\'{o}noma de Madrid.}
}
\date{} 

\maketitle

\begin{abstract}
We consider a model of fractional diffusion involving the natural nonlocal version of the $p$-Laplacian operator. We study the Dirichlet problem posed in a bounded domain $\Omega$ of $\ren$ with zero data outside of $\Omega$, for which the existence  and uniqueness of strong nonnegative  solutions is proved, and a number of quantitative properties are established. A main objective is proving the existence of a special separate variable solution $U(x,t)=t^{-1/(p-2)}F(x)$, called the friendly giant, which produces a universal upper bound  and explains the large-time behaviour of all nontrivial  nonnegative solutions in a sharp way. Moreover, the spatial profile $F$ of this solution solves an interesting nonlocal elliptic problem. We also prove everywhere positivity of nonnegative solutions with any nontrivial data, a property that separates this equation from the standard $p$-Laplacian equation.
\end{abstract}

\

\small
\setcounter{tocdepth}{1}
\begin{spacing}{0.0}
\tableofcontents
\end{spacing}
\normalsize
\newpage

\section{Introduction}  \label{sec.intro}

In this paper we consider a model of fractional diffusion involving  nonlocal operators of the $p$-Laplacian type, and we describe the main properties of the solutions of the evolution process
posed in a bounded domain $\Omega$ of $\ren$ with zero Dirichlet conditions outside of $\Omega$.

Our main contributions are the following: after posing the problem in the framework of semigroups generated by maximal monotone operators, showing existence and uniqueness of so-called strong solutions, and establishing the basic quantitative properties, we  prove the existence of a special separate-variable solution $U(x,t)$ of the form $U(x,t)=t^{-1/(p-2)}F(x) $, that we call the {\sl Friendly Giant.} Special solutions have played a great role in the development of very many areas in nonlinear analysis and applied mathematics. This happens in particular in the area of diffusion and heat transport, where there are some famous instances like the Gaussian kernel, the Barenblatt solutions, different traveling wave solutions, and other interesting examples with a lesser influence.

Our special solution plays a similarly important role for the problem we are addressing: it  produces a universal upper bound for the whole class of nonnegative solutions of the homogeneous Dirichlet problem. This means that the whole class of solutions is uniformly bounded at any given time $t>0$, and the universal bound decays to zero as time grows. Actually, the accuracy improves as $t\to\infty$ : $U(x,t)$ gives the first approximation the asymptotic behaviour of the all nonnegative solutions. Finally, we show that the spatial profile $F(x)$ of this solution solves an interesting nonlocal elliptic problem.

Let us be more precise: our aim is to study the equation
\begin{equation}\label{eq.1}
\partial_t u(x,t)+\int_{\ren} \frac{\Phi(u(y,t)-u(x,t))}{|x-y|^{N+sp}}\,dy=0,
\end{equation}
where $x\in\Omega\subset\ren$, $N\ge 1$, $\Phi(z)=c|z|^{p-2}z$, $1<p<\infty$, and $0<s<1$. The constant $c>0$ is unimportant.  The integral operator may be called the $(p,s)$-Laplacian operator, or more precisely the $s$-fractional $p$-Laplacian, our notation is ${\cal L}_{p,s}$. In this paper we cover the range $2<p<\infty$. Note that for $p=2$ we obtain the standard $s$-Laplacian heat equation, $u_t+(-\Delta)^{s} u=0$, about which much is known; on the other hand, it is proved that in the limit $s\to1$ with $p\ne 2$, we get the well-known $p$-Laplacian evolution equation $\partial_t u=\Delta_p(u)$, after inserting a normalizing constant. Due to this, the  notation $(-\Delta)^s_p$ is used in some papers for ${\cal L}_{p,s}$.

We recall that the standard $p$-Laplacian with $p>2$ is a well-known example of nonlinear diffusion of degenerate type, it is important in the study of free boundaries and as a benchmark for more complex nonlinear models, cf. the monograph \cite{DiBenBook, Lind}. The fractional version we are going to study has attracted attention only recently with the current focus on the analysis of nonlocal diffusion models of fractional Laplacian type, see the surveys \cite{CaffAbel, VazAbel}. References to previous work on equation \eqref{eq.1}  and its stationary versions are \cite{BCF, CafChanVas, ChLM, LL}. There is much current activity on this topic, see \cite{CKuusPal2015, IannMS2014, KMSire, MRT, Puhst, Ros15}.

\medskip

We will consider the equation in a bounded domain $\Omega\subset \ren$ with initial data
\begin{equation}\label{eq.ic}
u(x,0)=u_0(x), \qquad x\in \Omega,
\end{equation}
where $u_0$ is a nonnegative and integrable function. Moreover, we impose the homogeneous Dirichlet boundary condition that in the fractional Laplacian setting takes the form
\begin{equation}\label{eq.bc}
 u(x,t)=0 \qquad \mbox{ for all} \  x\in \ren, x\not\in \Omega,\quad \mbox{ and all $t>0$}.
\end{equation}
When we apply operator ${\cal L}_{p,s}$ on the set of functions that vanish outside of $\Omega$ we may be more precise and use the notation ${\cal L}_{p,s,\Omega}$, but this longer notation will not be used since there is no fear of confusion.
We call Problem DP the problem of solving  \eqref{eq.1} with conditions \eqref{eq.ic} and \eqref{eq.bc}. More precisely, we should write DP$({\cal L}_{p,s}; \Omega)$.

\medskip

\noindent {\bf Outline of the paper}

\noindent In Section \ref{sec.prel} we  gather the necessary preliminaries and notations, we analyze the differential operator ${\cal L}_{p,s}$, and we state the basic existence and uniqueness results that follow from standard theories once the problem is correctly posed. We also derive the basic properties of the  evolution semigroup corresponding to Problem DP.

In Section \ref{sec.bbdness} we state and prove  some of our main results: the existence of the special separate variable solution $U(x,t)$ that we call the Friendly Giant, and its role as universal upper bound, Theorems \ref{thm.1} and \ref{thm.2}. A historical source reference for this type of result is the paper by Aronson-Peletier \cite{AP81} on the Porous Medium Equation, 1981. The work was continued in \cite{DK88} where the curious name originated. A complete account of the results, for positive and sign-changing solutions was given in \cite{JLVmonats}, see also \cite{VazBook}.

The special profile $F(x)$ of our  solution $U(x,t)$ is the unique positive solution (groundstate) of an interesting nonlocal elliptic problem, see formula \eqref{form.3.2} in Theorem \ref{thm.1}.

Section \ref{sec.asbeh} deals with the asymptotic behaviour. We establish the correct rate $u(x,t)=O(t^{-1/(p-2)})$ and the sharp form $u(x,t)\,t^{1/(p-2)}\to F(x)$ in Theorem \ref{thm.3}. Consequently, we prove that the asymptotic behaviour does not depend at all on the solution (it is universal), as long as $u_0$ is nonnegative, nontrivial and integrable.

In Section \ref{sec.posit} we investigate whether the solutions have the finite propagation property of the standard $p$-Laplacian or, on the contrary, the nonlocal effect produces positivity at all points and for all times. The latter option is concluded: long distance effects win and nonnegative solutions are positive everywhere.

Finally, in Section \ref{sec.asbeh2} we give an improved version of the asymptotic convergence theorem, with sharp rate of convergence, under a boundary regularity assumption similar to the one established in \cite{RosSerra14, FRRO} for $p=2$. We state this assumption as an important conjecture.

The paper  concludes with a long list of comments and open problems.

\medskip

\noindent {\bf Fractional Laplacian operators, motivation and related works}

\noindent Let us  recall that fractional Laplacian operators and related nonlocal operators of different types have attracted considerable attention in recent times because of the interest in the applications and also because of the rich mathematical theory that has arisen. Indeed, interest was fueled  by a number of applications in various fields like continuum mechanics, stochastic processes of L\'evy type, phase transitions, population dynamics, image processing, finance, and so on. The original models were mostly linear, indeed the linear fractional heat equation has a long tradition in Probability since fractional Laplacians are  the infinitesimal generators of stable L\'{e}vy processes \cite{Applebaum, Bertoin, CKS2010}. However, a decade ago intense work began on nonlinear models, featuring  both stationary and  evolution models. Among the first, we mention the work on nonlocal obstacle models, \cite{CSS, Caffarelli-Silvestre}, on semilinear elliptic problems of nonlocal type \cite{CabSire, ServVald, ServVald2}, on phase transitions  \cite{SireVald}, and the work on nonlocal perimeters \cite{CaffRoqSav, CaffVald}.

Evolution problems of nonlocal diffusion type originated from material science and statistical mechanics, and led to much work on fractional porous medium models, see for instance \cite{BKM, CV1, CV2, CSV, pqrv1, pqrv2, BIK, BoVa14, STV14, VaVol14, CHSV15}, the work was summarized in \cite{VazAbel, VazRev2}.
 The Dirichlet Problem for  a version of the fractional Porous Medium  in a bounded domain has been recently treated by the author in collaboration with Bonforte and Sire, \cite{BSV, BoVa15}.
The close parallel that has been found between the theories of the standard Porous Medium Equation and the $p$-Laplacian equation, accredited in publications like \cite{DiBenBook, JLVSmoothing, IagSanVaz}, and the progress achieved in the study of fractional porous medium  models makes it natural to investigate the fractional $p$-Laplacian equation \eqref{eq.1}.

We will be inspired in the techniques developed to treat the Dirichlet Problem for the standard Porous Medium Equation, as explained in detail in \cite{JLVmonats}, where references to related work are found. The main ideas involved are rescaling, existence of special solutions, a priori estimates, and monotonicity arguments. The rescaled orbits converge to stationary states which solve a nonlinear elliptic problem. Such general techniques were then applied to the evolution $p$-Laplacian equation, and then the fractional Porous Medium.

The existence of a certain special solution is crucial in our analysis. Indeed, almost no special solution is known for  equation \eqref {eq.1}. In this respect we should mention Lindgren-Lindqvist's \cite{LL} who find the first eigenvalue and eigenfunction of the operator ${\cal L}_{p,s}$ on a bounded domain and describe its properties, in particular the limit as $p\to\infty$.

Finally, the unpublished paper \cite{MRT} contains the derivation of the basic existence and uniqueness results for different initial and boundary  value problems in the language of semi-group theory. Together with \cite{IshiiN, LL}, it allows us to shorten a bit our preliminaries section. In fact, the basic theory falls quite well into the scope of classical nonlinear semigroup theories for monotone or accretive operators, \cite{BrBk73, CL71}.

\medskip

\noindent {\sc Reminder and notations.} In the sequel $p>2$ and $0<s<1$ will be fixed. Only nonnegative solutions will be considered. The space-time domain of the solutions is denoted by $Q=\Omega\times (0,\infty)$. Also, $\overline{\Omega}$ denotes the closure of $\Omega$, and $d(x)=d(x,\partial \Omega)$ denotes the distance from a point $x\in\Omega$ to the boundary $\partial \Omega$. The notation $f_+$ denotes the positive part of a function, i.\,e.,  $\max\{f(x),0\}$. We write both $\partial_t u$ and $u_t$ for the partial derivative in time, either for clarity or brevity.

\section{Preliminaries}\label{sec.prel}

$\bullet$ Given numbers $p\in (1,\infty)$ and $s\in (0,1)$ we define the Gagliardo functional for smooth integrable functions defined in $\ren$ as
\begin{equation}\label{Jsp1}
J_{p,s}(u)=\frac1{p}\int_{\ren }\int_{\ren } \frac{|u(x)-u(y)|^p}{|x-y|^{N+sp}}\,dxdy\,.
\end{equation}
By extension it defines a norm in the fractional Sobolev space (also called  Gagliardo or Slobodetski space) $W^{s,p}(\ren)$. Indeed, the norm is given by
\begin{equation}\label{Jsp.norm}
\|u\|_{W^{s,p}(\ren)}^p =J_{p,s}(u)+\int_{\ren} |u(x)|^p\,dx\,.
\end{equation}
We refer to \cite{NPV} for information about these spaces and their properties. It is known that as $s\to 1$ the space $W^{s,p}$ becomes $W^{1,p}$, see \cite{BBM02}.

Functional $J_{p,s}$ is a convex, lower semi-continuous and proper funcional, and it has an associated Euler-Lagrange operator, given in its weak form by the expression
\begin{equation}\label{Lsp.weak}
  \int_{\ren }\int_{\ren } \frac{  |u(x)-u(y)|^{p-2}(u(x)-u(y))\, (\phi(x)-\phi(y))}{|x-y|^{N+sp}}\,dxdy
\end{equation}
where $\phi$ is any smooth variation in $W^{p,s}(\ren)$. This defines the operator ${\cal L}_{p,s}$ as
\begin{equation} \label{Lsp1}
{\cal L}_{p,s}(u)(x) =2 \, p.v. \int_{\ren }\frac{|u(x)-u(y)|^{p-2}(u(x)-u(y))}{|x-y|^{N+sp}}\,dy\,,
\end{equation}
where ``p.v.''  stands for the principal value of the integral. Thus defined, it is a positive operator, that corresponds for $p=2$ to the standard definition of
$(-\Delta)^{s}$ (up to a constant that we do not take into account). Caveat: we can find in the literature the definition with the opposite sign, to look like the usual Laplacian, but we will stick here to this definition as a positive operator, following the tradition in fractional Laplacian operators.

\medskip

\noindent $\bullet$ When restricted to the Hilbert space $L^2(\ren)$, ${\cal L}_{p,s}$  becomes a maximal monotone operator, as the sub-differential of functional $J_{p,s}$. It has a dense domain
$$
D({\cal L}_{p,s})=\{u\in L^2(\ren): J_{p,s}(u)<\infty, \ {\cal L}_{p,s}u\in L^2(\ren)\}\subset L^2(\ren).
$$
This means a number of properties for the evolution equation $\partial_tu+ {\cal L}_{p,s}u=0$. According to the general theory in  monograph \cite{BrBk73}, the usual mild solutions of the evolution equation provided by the semigroup theory are indeed strong solutions, which  means that for every $u_0\in L^2(\ren)$ the evolution orbit $u(t): C([0,\infty):L^2(\ren))$ is  differentiable for a.e. $t>0$ and $du/dt= -{\cal L}_{p,s}u\in L^2(\ren)$.  Besides,
    \begin{equation}
    \|u_t\|_{L^2(\ren)}\le \frac{C}{t}\|u_0\|_{L^2(\ren)}.
    \end{equation}
Moreover, a  number of properties of ${\cal L}_{p,s}$ are similar to the standard $p$-Laplacian operator, $-\Delta_p$. The convergence as $s\to 1$ of the solutions of ${\cal L}_{p,s}=f$ to the solutions of $-\Delta_p u=f$ is proved in \cite{IshiiN}. As a token of that the similarity, it is shown that  ${\cal L}_{p,s}$ is contractive in all $L^q$-norms, $1\le q\le \infty$, and this makes it into an $m$-completely accretive  operator in the sense of \cite{BCr91}, see proofs in \cite{MRT}. Consequently, the map  $S_t:u_0\mapsto u(t)$ is an $L^q$-contraction in $L^q(\ren)$ for every $1\le q\le \infty$ (and every $t>0$). We can also use the Implicit Time Discretization Scheme to generate the evolution solution by solving the elliptic problems
$$
h {\cal L}_{p,s}u+ u =f
$$
with $h>0$ a small constant (the time step), and applying the Crandall-Liggett iterative scheme  \cite{CL71}.

\medskip

\noindent $\bullet$ Let us now consider the problem posed on a bounded domain $\Omega\subset \ren$. For simplicity we will assume that $\Omega$ has a smooth boundary but this is condition can be relaxed. In order to take into account the homogeneous Dirichlet boundary conditions in the fractional setting, we restrict the operator to the set of functions $W_0^{s,p}(\Omega)$, defined as the closure of $C^\infty_c(\Omega)$\footnote{These functions are naturally extended by zero outside of $\Omega$.} with respect to the norm induced by $J_{p,s}$ according to   \eqref{Jsp.norm}. The formula for the sub-differential is the same, only the functional space changes. Indeed, formula \eqref{Lsp1} still applies as the defining formula with $x\in \Omega$, but integration in $y$ is  performed all over $\ren$. This means that we can feel the influence of the value $u(y)=0$ imposed for $y\not\in\Omega$ on the diffusion process inside $\Omega$, since the formula involves the differences $u(x)-u(y)$.

We thus finally obtain a maximal monotone operator defined in a dense domain in $L^2(\Omega)$ that is also $m$-completely accretive, see more in \cite{MRT}. We should denote this bounded domain operator as ${\cal L}_{p,s,\Omega}$ to avoid confusions with the operator defined for all $x\in\ren$, but we feel that this extra care is not needed, so ${\cal L}_{p,s}$ will be used in the sequel when we work on a bounded domain. In any case, there are many ways in which the difference between the action of ${\cal L}_{p,s}$ in the whole space and on a bounded domain is reflected, a very clear case is made in the first comment of Section \ref{sec.comm} concerning the lack of mass conservation for problem DP, for which we obtain a simple formula.

\medskip

$\bullet$ We get the following existence and uniqueness result.

\begin{thm} For every $u_0\in L^2(\Omega)$ there is a unique strong solution of Problem DP such that
(i) $u\in C([0,\infty):L^2(\Omega))$, (ii) for every $t>0$ ${\cal L}_{p,s}u(t)\in L^2(\Omega)$, and (iii) equation \eqref{eq.1} is satisfied for a.e. $x\in\Omega$ for every $t>0$. Moreover, for every two solutions $u_1, u_2$ with initial data $u_{01}, u_{02}$  we have
\begin{equation}
\|(u_1(t)-u_2(t))_+\|_1\le    \|(u_{01}-u_{02})_+\|_1\,, \quad\|u_1(t)-u_2(t)\|_q\le    \|u_{01}-u_{02}\|_q
    \end{equation}
for every $1\le q\le \infty$ and $t>0$.
 \end{thm}

Compare with Theorem  2.5 of \cite{MRT}. The theory can be applied to solutions with any sign. In the sequel we will deal with nonnegative solutions, $u(x,t), u_0(x)\ge 0$. We also point out that some references introduce the concept of viscosity solutions, like \cite{ChLM, IshiiN, LL, MRT}, but we will not need this concept here.

The class of initial data can be extended to $L^1(\Omega)$ using the $L^1$-contractivity. The resulting class of solutions is the same for $t>0$ since we will prove that all solutions are uniformly bounded for $t\ge \tau>0$, hence $u(t)\in L^\infty(\Omega)\subset L^2(\Omega).$

\medskip

 Here are some properties of the solutions that will be used below.

$\bullet$ Scaling property: a simple inspection shows that  ${\cal L}_{p,s}(A f(B x))= A^{p-1}B^{ps}({\cal L}_{p,s}f)(B x)$ for any $A,B>0$. When applied to the strong solutions of the evolution problem, it means that when $u(x,t)$ is a solution defined in the domain $\Omega$, then
    \begin{equation}
\widehat u(x,t)= A u(Bx,Ct)
\end{equation}
is a solution on $\Omega'=\{x: Bx\in \Omega\}$ on the condition that $C=A^{p-2}B^{ps}$. If we want to keep the domain fixed, then we need $B=1$. Notice that we also have ${\cal L}_{p,s}(- f( x))= -({\cal L}_{p,s}f)(x)$, so that $-u(x,t)$ is also a solution if $u$ is.

\medskip

$\bullet$ The collection of nonnegative solutions has the property of almost monotonicity in time
\begin{equation}\label{time.mon}
\partial_t u\ge -\frac{u}{(p-2)t}\,,
\end{equation}
valid a.e. in $\Omega$ for every $t>0$. The proof is based on the B\'enilan-Crandall homogeneity argument \cite{BCr}. Let us recall it briefly: scaling says that if $u_1(x,t)$ is a solution, also $u_\lambda(x,t)=\lambda u(x,\lambda^{p-2}t)$ is a solution. But if $\lambda>1$ we have $u_\lambda(x,0)\ge u_1(x,0)\ge 0$, hence by the maximum principle we get $u_\lambda(x,t)\ge u_1(x,t)$. Differentiate in $\lambda$ for $\lambda=1$ to get $(p-2)tu_t+u \ge 0$. A related homogeneity argument proves that
\begin{equation}\label{time.mon2}
\|\partial_t u\|_q\le \frac{2}{(p-2)t}\|u_0\|_q \quad \mbox{ for every } \ q\ge 1\,.
\end{equation}
 We ask the reader to complete this proof, see a similar argument in \cite[pg. 185]{VazBook}.

\medskip

$\bullet$ When $\Omega$ is a ball, say $B_R(0)$, and $u_0$ is radially symmetric, so is the solution $u(x,t)$ as a function of $x$. This is an immediate consequence of uniqueness and the rotation invariance of the operator.

\medskip

$\bullet$ Comparison of solutions in different domains applies: if $u_i$, $i=1,2$ are two solutions in domains $\Omega_i$ such that $\Omega_1\subset\Omega_2$ and we also have $u_1(x,0)\le u_2(x,0)$ for all
$x\in \Omega_1$, then $u_1(x,t)\le u_2(x,t)$ for all $x\in\Omega_1$ and $t>0$.

\medskip

$\bullet$ A usual property of nonlinear evolution processes of diffusion type is the so-called smoothing effect, whereby general initial data in a Lebesgue space $L^p(\Omega)$ will produce bounded solutions for positive times, with an $L^\infty(\Omega)$ that depends on the norm of the initial data. In the present case, much more will be proved, the bound will not depend on the initial data. This is what we study next.

\medskip

$\bullet$ Before going into that, let us state a result hat is not evident at all.

\begin{thm}\label{thm.pos} Let $u$ be the solution of problem DP corresponding to nontrivial initial data $u_0\in L^1(\Omega)$, $u_0\ge 0$. Then $u(x,t)$ is  positive in the whole domain $Q=\Omega\times(0,\infty)$, and uniformly positive on compact subsets.
\end{thm}

The proof of positivity will be postponed to Section \ref{sec.posit}. This result eliminates the possible existence of free boundaries, which were a typical feature of the standard $p$-Laplacian equation. For the fractional porous medium equation, a similar infinite propagation result holds in the model studied in \cite{pqrv1, pqrv2}, but free boundaries appear in the model  studied in \cite{CV1, CV2}.

\medskip

$\bullet$ The question of regularity of the stationary problems ${\cal L}_{p,s}(u)=f$ has been studied recently in \cite{CKuusPal2015, IannMS2014}. Theorem 1.1 of the latter paper states that when $f\in L^\infty(\Omega)$ then $u\in C^{\alpha}(\overline{\Omega})$ for some $\alpha>0$.
Due to the explained regularization effect of the evolution problem DP we have $u_t\in L^\infty(\Omega)$ uniformly for every $t\ge \tau>0$, hence the solutions of Problem DP are uniformly H\"older continuous w.r.t. the space variable for $t\ge \tau>0$.  $C^\alpha$ regularity is directly obtained by Sobolev embeddings when $sp>N$.


\section{The Friendly Giant. Boundedness of the solutions}\label{sec.bbdness}

 We will prove the following two results

\begin{thm} \label{thm.1} Let $p>2. $  There exists a special function $U(x,t)$ of the form
\begin{equation}\label{fg.sol}
U(x,t)=t^{-1/(p-2)}F(x)
\end{equation}
that is a strong solution of Problem \eqref{eq.1}-\eqref{eq.bc} for every $t>0$. The profile $F$ is positive in $\Omega$ and vanishes on $\partial \Omega$. It is a multiple of the unique positive solution of the stationary problem:
\begin{equation}\label{form.3.2}
 f\in W^{s,p}(\ren), \qquad {\cal L}_{p,s} f=f(x) \quad \mbox{in } \ \Omega, \qquad f(x)=0 \quad \mbox{ for } x\in \ren \setminus \Omega\,.
\end{equation}
 Moreover, $F$ is  $C^{\alpha}$-H\"older continuous  for some $\alpha>0$, and it behaves near the boundary like $O(d(x)^s)$.
\end{thm}

Note that $F$ depends on $\Omega$, $p$ and $s$. As a side remark, the special solution \eqref{fg.sol} satisfies $\partial_t U=-U/((p-2)t)$, which means that the a priori estimate \eqref{time.mon} is indeed optimal (with optimal constant).

This special solution  will be called the Friendly Giant as we mentioned above.  Strictly speaking, $U(x,t)$  is not an admissible solution of Problem DP since $U(x,0)=+\infty$ for every $x\in \Omega$. The easy remedy is to apply a small time delay and consider the solutions $U(x,t+h)$ for $h>0$. On the other, it is precisely the unbounded initial data what makes it an interesting object in the theory. Thus, we will prove that it is a universal bound for all nonnegative solutions of the Dirichlet problem, as the following result shows.

\begin{thm}\label{thm.2} For every  nonnegative solution $u(x,t)$ with initial data in some $L^p(\Omega)$, $1\le p\le \infty$, we have
\begin{equation}
u(x,t)\le U(x,t):=t^{-1/(p-2)}F(x)\,.
\end{equation}
\end{thm}

\subsection {New variables} In order to prove the theorems, we start by making the change of variables
\begin{equation}\label{form.chov}
v(x,\tau)=(a+t)^{1/(p-2)}u(x,t), \quad\tau=\frac1{p-2}\log(a+t)
    \end{equation}
 with a constant $a\ge 0$. Making the calculations and using the scaling property of ${\cal L}_{p,s}$, the new equation is then
\begin{equation}\label{v.eq}
\partial_\tau v+ (p-2) {\cal L}_{p,s} (v)=v\,.
\end{equation}
This is called the rescaled equation. If $a=0$ we have a small problem with the initial data since $t=0$ implies $\tau=-\infty$, so that the solution $v(x,\tau)$ is an eternal solution, defined for $-\infty<\tau<\infty$. But on the other hand, it has a clear advantage: the property of time monotonicity \eqref{time.mon} is equivalent to saying  that \ $\partial_\tau v\ge 0$ in the sense of distributions (at least). This is an important tool at our disposal.

\subsection {A uniform bound}

\begin{lem} \label{Wbound} There is a bounded function $W(x)>0$ such that for every solution of Problem DP we have
\begin{equation}
u(x,t)\le t^{-1/(p-2)}W(x)\,.
\end{equation}
\end{lem}

\noindent {\sl Proof. } It is a question of finding a suitable candidate. Here the candidate $W(x)$  is assumed to have a top  mesa of height 1 (i.e., a maximum value $W(x)=1$) in the ball of radius $2R$ (which is going to be large), $W$ is radially symmetric and decreasing for $|x|>2R$, and decreases to zero as $|x|\to\infty$ in a way to be chosen. The use of the integral formula \eqref{Lsp1} shows that $\Phi(x)={\cal L}_{p,s}W(x)$ is positive and smooth for $|x|<2R$ (i.\,e., at all points of maximum), and it is smooth outside. Hence, for $|x| < R$ we have ${\cal L}_{p,s}W(x)\ge C>0$. Since $W\le 1$ we get
$$
{\cal L}_{p,s}W(x)\ge C\,W(x) \quad \mbox{for} \ |x|\le R.
$$
Putting $\lambda W$ instead of $W$ we can fix the  constant $C$ to any given value, say $C=1$.

\noindent {\bf Remark.} We do not know if ${\cal L}_{p,s}W(x)$ is positive for $|x|>2R$, it depends on the tail; for compact support it is negative with decay $O(|x|^{-N+sp})$. This does not look so good, but such delicate information is not needed in our next argument, that goes as follows:

\medskip

\noindent $\bullet$ We perform the comparison step. We assume that $\Omega\subset B_R(0)$. We write
$$
\partial_\tau (v-W)+ {\cal L}_{p,s} (v)-{\cal L}_{p,s} (W)= v-{\cal L}_{p,s}(W)  \quad \mbox{ for} \ x\in\Omega.
$$
We multiply this equation by $p(v-W)$ where $p\in C^1(\RR)$ with  $p(s)=0$ for all $s\le 0$, $p(s)=1$ for $s\ge a>0$, and $p'(s)> 0$ for $0<s<a$. If $j$ is the primitive of $p$ with $j(0)=0$, we then  get
$$
\frac{d}{d\tau }  \int_{\ren} j(v-W)\,dx+ \int_{\ren} ({\cal L}_{p,s} (v)-{\cal L}_{p,s} (W))p(v-W)dx
= \int_{\ren} (v-{\cal L}_{p,s}W)p(v-W)dx.
$$
Since  $v(x)=0$ for $x\not\in \Omega$ and $W(x)>0$ everywhere, we have $p(v(x)-W(x))=0$ for all $x\not\in\Omega$, so that the integrals can be taken indistinctly over $\Omega$ or over $\ren$.  The second term of this display can be written in two parts. The first part is
$$
 \int_{\ren} \int_{\ren} \frac{(|v(x)-v(y)|^{p-2}(v(x)-v(y)))(p(v(x)-W(x))-p(v(y)-W(y)))}
{|x-y|^{N+sp}}\,dxdy\,,$$
and similarly the second part using ${\cal L}_{p,s} (W)$ instead of ${\cal L}_{p,s} (v)$. Subtracting them, we see that the whole term is positive.  Finally, we replace the last integral by
$$
\int_{\Omega} (v-{\cal L}_{p,s}W)p(v-W)dx\le \int_{\Omega} (v-W)p(v-W)dx\,,
$$
using the pointwise inequality ${\cal L}_{p,s}W\ge W$ in $\Omega$. Putting all together, we get
$$
\frac{d}{d\tau }  \int_{\ren} j(v-W)\,dx\le \int_{\Omega} (v-W)p(v-W)dx\,.
$$
We may let $p(s)$ tend the Heaviside function so that $j(s)\to s_+$. Putting now $X(\tau):=\int  (v-W)_+\,dx $, the limit of the last inequality gives
$$
X(\tau)\ge 0, \quad \frac{d X}{d\tau }\le X(\tau).
$$
\normalcolor If the initial data are zero, then the solution is always zero. But this means that $v(x,\tau)\le W(x)$ for all $x\in \Omega$ and all $\tau>0$.  This is equivalent to saying that $t^{-1/(p-2)}W(x)$ is an upper bound for the corresponding solution of Problem DP.

\medskip

\noindent $\bullet$ Can we ensure that $X(\tau)=0$ at some given time? It is not difficult when the DP starts with a function  $u_0$ is bounded and has compact support in $\Omega$. We proceed as follows we want to compare  $u(x,t)$ with $(t+\ve)^{-1/(p-2)}W(x)$ before doing the change of variables (i.\,e., the rescaling). The comparison is true  at $t=0$ if $\ve$ is small enough since $W$ is uniformly positive in $\Omega$. We now pass to the rescaled formulation with $a=\ve$ so that $t=0$ corresponds to some $\tau_0>-\infty$.
Since $X(\tau_0)=0$, we can apply the previous result to conclude that  $X(\tau)=0$ for all $\tau>\tau_0$. This means that
$$
u(x,t)\le (t+\ve)^{-1/(p-2)}W(x)\le t^{-1/(p-2)}W(x)\,,
$$
which ends the proof.

\noindent $\bullet$ We want to prove that the upper bound holds for any solution of Problem DP with data $u_0\in L^1(\Omega)$. In this case we approximate $u_0$ by compactly supported and bounded functions like in the previous step, get the uniform estimate for the approximations, and pass to the limit using the contraction property of the semigroup. The proof is concluded. \qed

\subsection{Alternative method using eigenfunctions}
In order to prove the upper bound, we may use as starting barrier the eigenfunction results that have been obtained by Lindgren and Lindqvist in \cite{LL}. They consider the minimization of the fractional Rayleigh quotient
$$
J_{p,s}(\phi):=\frac{\displaystyle\int_{\ren}\int_{\ren}\frac{|\phi(y)-\phi(x)|^p}{|x-y|^{N+sp}}\,dxdy}
{\displaystyle\int_{\ren} |\phi(x)|^p\,dx}
$$
among all functions $\phi$ in the class $C^\infty_0(\Omega)$, $\phi\not\equiv 0$. They assume that  $p\ge 2$ and $0<s<1$. They also write $\alpha p$ instead of $N+sp$, so when $0<s<1$ we get $N<  \alpha p < N+p$ (this change of notation is not important here). They call the infimum of $J_{p,s}(\Omega)$ the first eigenvalue, $\lambda_1$. It is positive and depends on $p,s$ and $\Omega$. Actually, it increases as we shrink $\Omega$ by a scaling argument. We are interested in  Theorem 5, Section 3,  of \cite{LL}, that says that:

\noindent {\bf Existence of the first eigenfunction.} {\sl There exists a non-negative minimizer $f\in W^{s,p}(\Omega)$ with $f\not\equiv 0$ in $\Omega$ and $f=0$ in $\ren\setminus \Omega$. It satisfies the Euler-Lagrange equation
\begin{equation}\label{eigen}
{\cal L}_{p,s} f=\lambda_1 |f|^{p-2}f
\end{equation}
in the weak sense with test functions $\phi\in C^\infty_0(\Omega)$.
 If $sp>N$, the minimizer is in $C^{0,\beta}(\ren)$ with $\beta = s-{N/p}.$}

\noindent $\bullet$   We use the previous result to construct a supersolution of the eigenvalue problem
$$
{\cal L}_{p,s} (H)=c\,H
$$
inside a ball $B_R$ by using the first eigenfunction $f$ of the minimization problem posed in a larger domain, like the ball $B_{2R}. $ This will be a way to obtain the upper for $v$ with a more explicit construction.

In this way we may retrace our steps and prove the existence of the universal upper bound if $p$ is large, more precisely $sp>N$. \cite{LL} does in that case a viscosity theorem and proves uniqueness and positivity of the first eigenfunction.  But for smaller $p$ we have the regularity and positivity results that have been proven recently.

\subsection{The bounded limit profile $F(x)$}

Let us go back to our rescaled solution $v(x,\tau)$.  We have monotonicity time of the whole orbit for $-\infty<t<\infty$. We also have a uniform a priori bound $v(x,\tau)\le W(x)$ valid for all $x\in \Omega$ and all $\tau\in\RR$. Therefore, we are allowed to pass to the limit $t\to\infty$ and get the expression
\begin{equation}
\lim_{t\to \infty} v(x,\tau)=F(x)
\end{equation}
with $F(x)\le W(x)$. Since we also have a uniform bound of the form on $\partial_\tau u$ and ${\cal L}_{p,s}u$, it is not difficult to prove that the equation \eqref{v.eq} becomes
\begin{equation}
{\cal L}_{p,s} (F)=\mu\,F, \quad \mu=\frac1{p-2}\,,
\end{equation}
at least in the weak sense. The factor $\mu$ is not essential. It can be eliminated by a simple rescaling that uses the different homogeneities of the right-hand and left-hand side; is we put $F=cF_1$ and
$c^{p-2}=\mu$, then ${\cal L}_{p,s} (F_1)=F_1$.

Finally, since Theorem  \ref{thm.pos} asserts that $u(x,t)$ is positive everywhere, so is $v(x,\tau)$. $F(x)$ is a monotone increasing limit of $v$, so $F(x)>0$ for every $x \in\Omega$. We can now apply the regularity results of \cite{CKuusPal2015, IannMS2014} to conclude that $F$ is \ $C^{\alpha}$ H\"older continuous for some $\alpha>0$ and that it has the stated behaviour at the boundary (see \cite[Theorem 4.4]{IannMS2014}).

\subsection{Universal bound and uniqueness of the profile}

By means of passage to the limit of a particular nontrivial orbit we have obtained a Friendly Giant of the form $U(x,t)=t^{-1/(p-2)}F(x)$ that is an upper bound for the whole orbit $u(x,t)$ with initial data $u_0(x,t)$ from which we started. But the argument of Lemma \ref{Wbound} applies to any other solution $u_1(x,\tau)$ with initial data $u_{01}(x)\ge 0$ using now $F$ as an upper barrier in the place of $W$. In this way we conclude that $t^{-1/(p-2)}F(x)$ is a universal upper bound for any orbit.

Following this argument we can prove the uniqueness of the strong nonnegative solution of problem \eqref{v.eq}. Suppose we have a second solution $\widehat F\ge 0$. We consider the function $\widehat u$ solution of \eqref{eq.1} with data $\widehat F$ at $t=0$.  By the previous argument we have
$$
\widehat u(x,t)=(t+1)^{-1/(p-2)}F_1(x)\le U(x,t)=t^{-1/(p-2)}F(x)
$$
passing to the limit we conclude that $F_1(x)\le F(x)$. The roles can be reversed, hence $F_1=F$
The proof  of Theorems \ref{thm.1} and  \ref{thm.2} is done, pending the proof of positivity of the solutions. \qed

\medskip

\noindent {\bf Remarks.} (1) We have proved that equation \eqref{v.eq} for $v(x,\tau)$ has just one positive stationary solution. However, it has many non-stationary ones, derived from it and following the formula
$$
V(x,\tau;a)=(1+\frac{a}{t})^{1/(p-2)}F(x)\,.
$$
All of them tend obviously to $F$ as $t\to\infty$.

\medskip

\noindent (2) If we expand the domain by means of a homothetical scaling from $\Omega$ to  $\lambda \Omega=\{\lambda x: \ x\in \Omega\}$, then we get the following formula for the corresponding profile
\begin{equation}
F_{\lambda\Omega}(x)=\lambda^{ps/(p-2)}\,F_\Omega(x/\lambda)\,.
\end{equation}

\section{Asymptotic behaviour}\label{sec.asbeh}

Here, we establish the asymptotic behaviour of all nonnegative solutions of the Dirichlet problem.

\begin{thm}  \label{thm.3} For every nonnegative solution $u(x,t)$ of the DP we have
\begin{equation}
\lim_{t\to \infty} t^{1/(p-2)}\| u(x,t)-U(x,t)\|_q= 0\,.
\end{equation}
for every $L^q$ norm with $1\le q \le  \infty$. In other words, $t^{1/(p-2)}\, u(x,t)$ tends to $F(x)$ as $t\to\infty$ uniformly in $\Omega$.
\end{thm}

\noindent {\sl Proof.}  The statement is equivalent to saying that $v(x,\tau)$ converges uniformly to $F(x).$ Now, the time monotonicity and the boundedness of $F$ allows us to use the monotone convergence theorem to conclude that $v(\cdot, \tau)$ converges to $F(\cdot)$ in $L^1(\Omega)$. This and the uniform boundedness imply the convergence in all $L^q(\Omega)$ with $q<\infty$.  The extra result of uniform H\"older regularity of $v(\cdot,\tau)$ implies convergence in $L^\infty(\Omega)$ by standard interpolation results between $L^1$ and $C^\alpha$.
\qed

\medskip

\noindent {\bf Remark.}  We have only found a related result for the asymptotic behaviour of the present Dirichlet problem, recently announced in \cite{MRT}, that provides a first rough estimate. Indeed, their  Theorem 2.6 says:

\noindent $\bullet$ {\sl Let $q \ge p$. Let $u(x,t)$ be the solution of the Dirichlet problem  for the initial datum $u_0 \in L^\infty(\Omega$), if $q > p$ and $u_0 \in L^2(\Omega)$ if $q = p$. Then, the $L^q$-norm of the solution goes to zero as $t\to\infty$ since we have the following estimate:}
\begin{equation}
\|u(t)\|^q_{L^q(\Omega)}\le\frac{C}{t}\|u_0\|_{L^\infty(\Omega)}^{q-p}\|u_0\|_{L^2(\Omega)}^2
\end{equation}

\section{Positivity. Reflection principle}\label{sec.posit}

The proof of Theorem \ref{thm.pos} has two steps. In the first we prove the result for the case where
$\Omega$ is a ball and $u_0$ is radially symmetric. In the second we obtain the whole result.

\medskip

\noindent $\bullet$ {\sl Radial monotonicity.} Before we prove positivity of the solutions when $\Omega$ is a ball is $u_0(x)$ is radially symmetric around the center of the ball and decreasing with the distance, we need to prove a monotonicity principle in the radial direction. We may assume that $\Omega=B_R(0)$. We know that in that case the solutions are radially symmetric. We want to prove that they are monotone in the radial direction. The comparison is done by an application of Aleksandrov's reflection principle. This has been performed in the case of the fractional porous medium in  \cite[Section 15]{VazBar}, and we only have to repeat the process followed there.

 A main reduction is that we only need to prove the reflection principle for the elliptic steps of the implicit time discretization, which amounts to prove the result for solutions of
\begin{equation}\label{ell.1}
{\cal L}_{p,s} u + u =f\,.
\end{equation}
posed in $\Omega$ with $f\in L^2(\Omega)$.  We take a hyperplane $H$ that divides $\ren$ into two half-spaces, and let  $\Omega_1$ and $\Omega_2$ the two pieces into which $H$ divides $\Omega$. Let us assume that the symmetry $\Pi$ with respect to $H$  that maps $\Omega_1$ into $\Omega_2$ (i.\,e., $H$ does not pass through the origin, which is contained in $\Omega_2$). The statement we need to prove is as follows.
\begin{thm} Let $u$ be the unique solution of \eqref{ell.1} with data $f\in L^2(\Omega)$, $f\ge 0$. In the above circumstances and under the further assumption that
\begin{equation}
f(x)\le f(\Pi(x)) \qquad \text{in} \quad \Omega_1
\end{equation}
we have $u(x)\le u(\Pi(x))$ in $\Omega_1.$
\end{thm}

The proof does not differ essentially from the one in  \cite{VazBar}. As a corollary, we conclude that
whenever $f$ is radially symmetric and nonincreasing in $|x|$, so it $u(x)$. Passing to a similar statement about the parabolic evolution is immediate.

\medskip

\noindent $\bullet$ {\bf I}. We can now prove positivity for the solutions of the parabolic problem when $\Omega$ is a ball if $u(x,t)$ is radially symmetric around the center of the ball and decreasing with the distance (at all times).

Suppose that $u(r,t)$ is a radially symmetric solution, $r=|x|$, and it is nonincreasing in $r$. If for some $t>0$ have $u(r_1,t)=0$ a point $r_1<R$ we will also have $u(r,t_1)=0$ for all $r_1<r<R$. By the time monotonicity we will also have $u(r,t)=0$ in the cylindrical region: \{$r_1<r<R$, $0<t<t_1$\}. Now, in that region $u_t=0$ while ${\cal L}_{p,s}u<0$, since the points are points of minimum of $u$, and $u$ is not identically zero. This is a contradiction with equation \eqref{eq.1}, hence $u$ cannot touch zero. We conclude that $u(r,t)>0$ everywhere in $Q$.

\medskip

\noindent $\bullet$ {\bf II. \sl Positivity for general solutions.} We take now a general domain and only assume on $u_0$ that is strictly positive in some ball, say $u_0(x)\ge \ve>0$ in  $B_R(x_0)\subset \Omega$. We consider the solution $u_1(x,t)$ with domain $B_R(x_0)$ and initial data $u_{01}$ smooth, radially symmetric around $x_0$, decreasing in $r=|x-x_0|$. By the previous item, $u_1(x,t)$ is positive
everywhere in $Q(R;x_0):=B_R(x_0)\times (0,\infty)$. A simple application of the maximum principle implies that $u(x,t)\ge u_1(x,t)$ in $Q(R;x_0)$. Notice that maximum radius $R$ we may take  the distance from $x_0$ to the boundary $\partial \Omega$.

In order to propagate this positivity everywhere in $Q=\Omega\times (0,\infty)$, we select a point
$x_1\in \Omega$ and a time $T>0$ and we argue as follows: if $|x_1-x_0|$ is less than the distance from $x_0$ to the boundary we are already done. If this is not the case there is a finite set of points
$y_1, y_2, \cdots, y_N=x_1$, such that $|y_{k+1}-y_{k}|< \frac1{2} d_k$, where $d_k= \mbox{dist}(y_k,\partial \Omega)$.
Picking times $t_0=0<t_1<t_2<\cdots<t_n=T$ and applying iteratively the previous argument, we conclude that $u(x, t_k)>\ve_k>0$ for $x\in B_{d_k/2}(y_k)$. In the end $u(x_1, T)>\ve_N>0$.

Finally, when $u_0$ is not known to be positive in a ball, we wait a time $\ve>0$ and apply the regularity results to show that $u(x,\ve)$ must be positive in a ball since it is continuous, nonnegative and nontrivial. The rest of the argument is the same for $t>\ve>0$. Since $\ve$ is arbitrary, we get positivity everywhere for all $t>0$. \qed

\section{Sharp boundary behaviour, asymptotics with rate}\label{sec.asbeh2}

In a remarkable paper \cite{RosSerra14} Ros-Ot\'on  and Serra studied the regularity up to the boundary of solutions to the Dirichlet problem for the fractional Laplacian, case $p=2$. They  consider the  solutions $u$ of the elliptic equation $(-\Delta)^s u=f$ with $0<ws<1$ and $f\in L^\infty(\Omega)$, and satisfying Dirichlet condition $u=0$ for $x\not\in \Omega$. They prove that $u\in C^s(\ren)$ and $u(x)/d(x)^s$ is $C^\alpha$  in $\Omega$ and up to the boundary, for some $\alpha (0,1)$, where $d(x) = dist(x; \partial\Omega)$. They do that by developing a fractional analog of the Krylov boundary Harnack method. This assertion must be strengthened with the assertion that the quotient $u(x)/d(x)^s$
is strictly positive, which is done in Lemma 3.2 of the paper, see also \cite{Ros15}. The results were extended to the fractional heat equation in \cite{FRRO}.

\noindent $\bullet$ We conjecture that this boundary behaviour result is true for $p>2$ with the same power, $u(x,t)\sim d(x)^s$. The upper bound for $F(x)/d(x)^s$ is true by the already mentioned upper boundary behaviour estimate,
\cite{IannMS2014}, and as a consequence of the uniform bound in terms of $U(x,t)$, the same upper bound
applies to any solution $u(x,t)$. But the sharp lower bound is missing at this moment, and it has important consequences.

\medskip

If this is accepted as a sensible hypothesis, we can prove the following version of the sharp asymptotic convergence with rate.

\begin{thm}  \label{thm.4} Let the conjecture hold and let  $u$ be nonnegative solution $u(x,t)$ of problem  DP with nontrivial data. Then there is a constant $T>0$ that depends on the data such that
\begin{equation}
 U(x,t+T)  \le  u(x,t)  \le U(x,t)
\end{equation}
holds in $\Omega$ if $t$ is large. This means that as $t\to\infty$
\begin{equation}
t^{1/(p-2)}u(x,t)= F(x)+ O(1/t).
\end{equation}
\end{thm}

\noindent {\bf Remark.} The example $u(x,t)=U(x,t+T)$ shows that the result is sharp.

\medskip

\noindent {\sl Proof.} We only need to establish the lower bound,  $U(x,t+T)  \le  u(x,t) $. By the maximum principle, this can be done at any given time, say $t=t_1$ and then it will hold for all times $t\ge t_1$. Now, at $t=t_1$ we have by the assumed regularity
$$
u(x,t_1)\ge C\,d(x)^s, \quad C=C(u,t_1)
$$
while away from the boundary $u(x,t_1)$ is uniformly positive by continuity and the positivity result of the previous section. We conclude that there is an $\ve>0$ such that $
u(x,t_1)\ge \ve\,F(x).$
On the other hand,  for $T$ large enough
$$
U(x,t)=(t_1+T)^{-1/(p-2)}F(x)\le \ve F(x)\,.
$$
The conclusion follows. \qed
\normalcolor

\section{Comments and open problems}\label{sec.comm}

We list here a number of  comments and questions that might interest the curious reader.

\noindent $\bullet$  The difference between the action of  operator ${\cal L }_{p,s} $ in $\ren$ and on a bounded domain can be seen when we compute the lack of conservation of mass for a strong nonnegative solution of problem DP  in a bounded $\Omega$. From the equation we have
$$
\frac{d}{dt}\int_\Omega u(x,t)\,dx=-\int_\Omega {\cal L }_{p,s}u(x,t)\,dx.
$$
After putting $1_\Omega(x)=\phi(x)$ and using the fact that $u(x,t)$ can be extended by zero if $x\not\in \Omega$, the last integral equals, a
$$
\int_\ren {\cal L }_{p,s}u(x,t)1_{\Omega}\,dx=\int_{\ren }\int_{\ren } \frac{  |u(x)-u(y)|^{p-2}(u(x)-u(y))\, (\phi(x)-\phi(y))}{|x-y|^{N+sp}}\,dxdy\,.
$$
We only need to take into account the case where $x\in\Omega$ and $y\not\in\Omega$ (and the converse), otherwise $\phi(x)-\phi(y)=0$.
Therefore, we get
$$
\frac{d}{dt}\int_\Omega u(x,t)\,dx=-2\int_{x\in\Omega} \int_{y\not\in\Omega}
\frac{|u(x)|^{p-2}u(x)}{|x-y|^{N+sp}}\,dxdy,
$$
which is never zero, and determines the loss of mass. In $\ren$ such an integral does not exist
and mass is conserved for all solutions for which the calculation is justified.

\noindent $\bullet$  The regularity of  strong solutions of Problem DP is a very ongoing interesting question. See the works  \cite{CafChanVas, Kass09} for non-degenerate problems of this type, where $C^\alpha$ regularity of bounded solutions is proved. For fractional porous medium equations it has been proved in  \cite{CSV, CVUral} and \cite{VPQR} for the different models. As we have noted above, sharp boundary regularity is really needed.

\noindent $\bullet$ The universal upper bound is true for signed solutions, that are naturally produced by the theory when $u_0$ is not necessarily negative. We get in particular the following version for a signed solution: $|u(x,t)|\le t^{-1/(p-2)}F(x)$. This translates into corresponding bounds for the large-time decay  of the form $|u(x,t)|=O(t^{-1/(p-2)})$, but the sharp versions of the asymptotic behaviour are lost since $F(x)$ is not representative, generally speaking. See the analysis for the standard PME and PLE in \cite{JLVmonats}.

\noindent $\bullet$ Theorem \ref{thm.2} states that the universal upper bound applies to all solutions with initial data in $L^1(\Omega)$. By passage to the limit it would also apply to solutions having as initial data a nonnegative Radon measure. We will not enter into the theory for such solutions, since more basic questions are still open.

\noindent $\bullet$ The universal upper bound and universal asymptotic behaviour are true for all $p\in(2,\infty)$, but they do not hold for the limit values of $p$. Thus, they cannot be true for the limit case $p=2$, since it is a linear equation and the size of the solutions must depend linearly on the size of the data. For the other limit, $p\to\infty$, an evolution theory has to be developed, see \cite{ChLM, LL} for the proposed elliptic operator. The reader will notice that the time factor $t^{-1/(p-2)}$ of the Friendly Giant goes to 1 in the limit, its decay disappears.

\noindent $\bullet$ Regarding the asymptotic behaviour in the linear case $p=2$, it is known that solutions of the Dirichlet Problem   decay exponentially, and more precisely, the large time behaviour of a solution takes the classical form $u(x,t)\sim c(u_0)\phi(x)e^{\lambda_1 t}$ where $\phi(x)>0$ is the first eigenfunction and $\lambda_1 >0$ is the first eigenvalue of ${\cal L}_{2,s}=(-\Delta)^s$. So the asymptotics is quite different.

\noindent $\bullet$  The large time  situation is even more different when $1<p<2$. Though we will not study that range in this paper since it deserves proper attention, let us just point out that the method of Section \ref{sec.bbdness} shows in a simple way that there a phenomenon of complete extinction in finite time, as the one already described for fast diffusion and fast $p$-Laplacian in the literature, \cite{JLVSmoothing}. A brief argument is as follows: we introduce a rescaled solution $v(x,\tau)$ by means of new formulas
\begin{equation}
v(x,\tau)=(T-t)^{-1/(2-p)}u(x,t), \quad  \tau=\frac1{2-p}\log (1/(T-t)),
\end{equation}
where $T>0$ is arbitrary and the ranges of time are $0<t<T$ and $\tau_0<\tau<+\infty$. This change leads to the equation for $v$:
\begin{equation}\label{v.eq2}
\partial_\tau v+ (2-p) {\cal L}_{p,s} (v)=v\,.
\end{equation}
Like in Section \ref{sec.bbdness}  we find that there is stationary positive supersolution $W(x)$ for this equation. If moreover the data or the solution are bounded, i.\,e., if  $v(x,\tau_1)\le W(x)$ for some $\tau_1$, we conclude that $v(x,\tau)\le W(x)$
for all $x\in\Omega$ and $\tau_1<\tau<  \infty$. Translating this for $u$ we get
$$
u(x,t)\le (T-t)^{1/(2-p)}W(x)\quad \mbox{for \ } x\in\Omega, \ t_1<t<T\,,
$$
which implies that $u(\cdot,y)$ converges uniformly to zero as $t\nearrow T$. We conclude that the orbit $u(\cdot,t)$ vanishes identically after a time $T(u)$ which is equal or less than the constant $T$
of the supersolution. Analyzing the behaviour for $t$ near the extinction time is an open problem.

\noindent $\bullet$  The problem with right-hand side, $\partial_t u+{\cal L}_{p,s} u=f$ is interesting. Existence and uniqueness of solutions is granted by the usual  theory.  See \cite{KMSire} for recent results on the elliptic problem.

\noindent $\bullet$ Some authors have considered operators like ${\cal L}_{p,s}$ with more general nonlinearities than our $\Phi(z)=|z|^{p-2}z$. Our results do not apply to them, at least the sharp ones.

\noindent $\bullet$ There are alternative definitions of the fractional Laplacian on bounded domains in the linear case $p=2$, that have been discussed recently in the literature, see e.\,g. \cite{BSV, MuNaz14, SVspct}. No such alternative definition seems to apply in our $p$-Laplacian case. This is an intriguing question.

\noindent $\bullet$ It would be nice to know more semi-explicit solutions of equation \eqref{eq.1}.

\noindent $\bullet$ Note that the evolution generated by Problem DP in $\Omega$, started with initial condition the eigenfunction \eqref{eigen} of \cite{LL}, produces in the rescaled $v$-formulation a natural connection between this profile and the friendly giant profile $F(x)$.

\noindent $\bullet$ The methods we have used do not allow us to settle the fine asymptotic behaviour of the solutions of the Cauchy problem posed on $\ren$. The idea is to construct the corresponding Barenblatt solutions, but such process is not easy.

\noindent $\bullet$ Can we say that our $p$-$s$ Laplacian equation with $p>2$ is a degenerate parabolic equation when $0<s<1$, as it is said of the same equation when $s=1$?

\vskip .8cm

\noindent{\sc Acknowledgements.} Work partially funded by Project MTM2011-24696 (Spain). The author is grateful to Y. Sire for very valuable information on current literature, and to some  other colleagues for helpful comments.

\



    \normalsize

\vskip 1cm

\noindent {\bf Keywords.} Nonlinear evolutions, $p$-Laplacian equation, fractional diffusion.

\medskip

\noindent{\bf Mathematics Subject Classification}. 35B45, 35B65,
35K55, 35K65.\\


\vskip .5cm

\noindent {\sc Address:}

\smallskip

\noindent {\sc Juan Luis V{\'a}zquez}\newline
Departamento de Matem\'{a}ticas, Universidad Aut\'{o}noma de Madrid, \\
28049 Madrid, Spain. \ e-mail: {\tt juanluis.vazquez@uam.es}

\vskip 1cm

\end{document}